\documentclass[journal]{IEEEtran}
\usepackage{amsmath}
\usepackage{amsfonts}

\newtheorem{theorem}{Theorem}

\newtheorem{remark}[theorem]{Remark}

\input{tcilatex}
\begin{document}

\title{Mathematical Representation \\
of Multitarget Systems}
\author{Ronald Mahler, Random\ Sets LLC, Eagan, MN, U.S.A.}
\maketitle

\begin{abstract}
This paper systematically compares two mathematical foundations for
multitarget tracking: \ labeled random finite sets (LRFS's) and trajectory
random finite sets (TRFS's).
\end{abstract}

\section{Introduction \label{A-Intro} \ }

As an engineering discipline, multitarget tracking dates back to the seminal
paper by Reid \cite{Reid-MHT} and earlier---see \cite{WileyEnc2015} for an
overview. \ The last 20 years of multitarget tracking research have differed
from the 20 years preceding them in that there has been an increasing
emphasis on theoretically rigorous statistical foundations that also
facilitate the development of practical multitarget tracking algorithms. \ 

Probably the most notable such algorithm is the generalized labeled
multi-Bernoulli (GLMB) filter, which arose from the random finite set (RFS) 
\cite{Mah-Artech} and labeled random finite set (LRFS) \cite%
{Vo-ISSNIP12-Conjugate}, \cite{VoVoTSPconjugate} paradigms. \ Its most
recent implementations---made possible by the application of advanced Gibbs
statistical sampling techniques to the LRFS multitarget posterior density \ $%
f_{k}(X|Z_{1:k})$ \ on labeled multitarget state-sets \ $X$---can
simultaneously track over one million targets in significant clutter in real
time using off-the-shelf computing equipment \cite{Beard2020}.

In the LRFS approach, the state of a multitarget system at time \ $t_{k}$ \
is modeled as a \textit{labeled }finite subset \ $X=\{(\mathbf{x}_{1},\ell
_{1}),...,(\mathbf{x}_{n},\ell _{n})\}\subseteq \mathfrak{X}_{0}\times 
\mathfrak{L}$ \ where: \ $\mathfrak{X}_{0}$ \ is the kinematic state space
and \ $\mathfrak{L}$\ \ is a countable set of \textquotedblleft track
labels\textquotedblright ; $\ X_{\mathfrak{L}}=\{\ell _{1},...,\ell _{n}\}$
\ is the set of labels of the elements of \ $X$; and $|X|=|X_{\mathfrak{L}}| 
$ \ where \ $|S|$ \ denotes the cardinality of a finite set \ $S$. \ That
is, \ $X$ \ is \textquotedblleft labeled\textquotedblright\ if the targets
in \ $X$ \ have unique labels. \ Given this, a time-evolving\ multitarget
population is represented as a time-sequence \ $X_{1},...,X_{k},...$ \ of
multitarget labeled state-sets at the measurement-collection times \ $%
t_{1},...,t_{k},...$ \ 

Since 2014, however, it has repeatedly been asserted that the LRFS framework
is seriously deficient because\ LRFS labels in \ $\mathfrak{L}$ \ have
\textquotedblleft ...no direct and unambiguous connection to the physical
phenomena under observation...\textquotedblright\ (\cite{Svensson-FUSION2014}%
, p. 2, column 2). \ Worse, \textquotedblleft ...adding labels to the state
artificially increases the uncertainty in the tracking
problem...\textquotedblright\ (\cite{Svensson-FUSION2014}, p. 3, top of
column 1).

The remedy that has been proposed is as follows: \ Replace the LRFS
framework with an allegedly \textit{label-free} \textquotedblleft trajectory
random finite set\textquotedblright\ (TRFS) framework. \ That is, replace
the time-sequence \ $X_{1},...,X_{k}$ \ with a finite \textquotedblleft set
of trajectories\textquotedblright\ (SoT) \ $\mathbf{T}=\{T_{1},...,T_{n}\}$
\ where each \textquotedblleft trajectory\textquotedblright\ \ $T_{k}$ \ is
a vector of the form \ \ 
\begin{equation}
T_{k}=(t_{k},\mathbf{x}^{1:i})\overset{_{\text{abbr.}}}{=}((t_{k},\mathbf{x}%
^{1}),(t_{k+1},\mathbf{x}^{2}),...,(t_{k+i-1},\mathbf{x}^{i}))
\end{equation}%
where \ $t_{k}$ \ and \ $t_{k+i-1}$ \ are the trajectory's beginning and end
times and where \ $\mathbf{x}^{1:i}:\mathbf{x}^{1},\mathbf{x}^{2},...,%
\mathbf{x}^{i}\in \mathfrak{X}_{0}$ \ are its respective kinematic states at
those times.

The purpose of this paper is to substantiate the following claims: \ 

\begin{enumerate}
\item These criticisms of the LRFS framework are mistaken;

\item the TRFS framework itself has no direct and unambiguous connection to
the physical phenomena under observation; and

\item it itself seriously increases uncertainty in the tracking problem. \ 
\end{enumerate}

The analysis that follows will begin with phenomenological and
epistemological basics---specifically, with the basic semiotic concepts of 
\textit{signified}\ (an entity) versus \textit{signifier}\ (a name for that
entity in some symbolic language). \ It will be shown, step-by-step, that
the LRFS framework arises inevitably from the phenomenological requirements
of the multitarget tracking application. \ It will then be shown that the
TRFS framework arises by stripping from the LRFS\ framework the signifiers
(labels) that signify \textit{actual physical phenomena}---i.e., the unique
identities that real-world targets actually possess. \ Stated more
forthrightly: \ 

\begin{enumerate}
\item[$\bullet $] Jettisoning the signifier (the label and the information
that it carries) effectively throws the baby (the unique target that it
signifies) out with the bathwater.
\end{enumerate}

The remainder of the paper is organized as follows: \ mathematical
representation of targets (Section \ref{A-TargRep}); mathematical
representation of target trajectories (Section \ref{A-TrajRep}); and
conclusions (Section \ref{A-Concl}).

In what follows, \ $\mathbb{N}=\{0,1,2,...\}$ \ denotes the natural numbers;
\ $\mathbb{N}^{+}=\{1,2,...\}$ \ the nonzero natural numbers; \ $\mathbb{R}%
^{+}$ \ the positive real numbers; and \ $\delta _{i,j}$ \ the Kronecker
delta.

\section{Mathematical Representation of Targets \label{A-TargRep}}

The section is organized as follows: \ mathematical representation of target
kinematics (Section \ref{A-TargRep-AA-Kine}); mathematical representation of
target identity (Section \ref{A-TargRep-AA-Ident}); and mathematical
representation of multitarget states (Section \ref{A-TargRep-AA-MultiTarg}).

\subsection{Mathematical Representation of Target Kinematics \label%
{A-TargRep-AA-Kine}}

A target is a macroscopic physical entity. \ This means in particular that
it has a position, which is typically mathematically idealized as a unique
point in three-dimensional Euclidean space. \ A point position must be
designated by some sort of symbolic identifier---a \textquotedblleft
signifier,\textquotedblright\ in the patois of semiotics. \ This signifier
is not---indeed, cannot be---unique. \ It typically consists of a triplet of
real numbers (continuous state variables) in some coordinate system
(Cartesian, cylindrical, spherical, etc.) centered at some origin point
(geocentric, etc.), each variable having a unit of measurement (metric,
English, etc.), and the value of each variable being expressed in terms of
some number-system base (decimal, binary, hexadecimal, etc.). \ 

Two points should be emphasized:

\begin{enumerate}
\item The arbitrariness of a position signifier does not negate the reality,
physicality, or uniqueness of the position that it signifies; and the same
is true of any other kinematic state variable. \ 

\item The existence of an infinitude of possible signifiers for a position
does not contravene the usual requirement that there be a one-to-one
correspondence between physical states and their mathematical
representations. \ For position (or any other kinematic state variable), it
is enough that this correspondence exists once the following have been
selected: \ number base, measurement units, coordinate system origin, and
coordinate system.
\end{enumerate}

\subsection{Mathematical Representation of Target Identity \label%
{A-TargRep-AA-Ident}}

As a macroscopic physical entity, a target has a discrete-valued state
variable---its unique identity---which must be assigned some appropriate
signifier. \ As with position, this signifier cannot be unique. \ For \
example, it can be (in the case of people) a social security number or a
name such as \textquotedblleft Bob Aardvark\textquotedblright\ or
\textquotedblleft Sue Zebra\textquotedblright ; or (in the case of aircraft)
a tail number or a manufacturers' serial number. \ 

The arbitrariness of this signifier does not preclude the reality,
physicality, or uniqueness of the target identity which it signifies. \ And
the infinitude of possible identity-signifiers does not violate the
one-to-one correspondence principle.\ \ Since the number of target
identities is finite but usually unknown, it is enough that such a
correspondence exists once a countable set of identity-signifiers has been
chosen.

When the actual identity of a physical object is unknown, it can be assigned
a \textit{tentative identifier}---for example, a \textquotedblleft track
label.\textquotedblright\ \ Neither the provisionality nor the infinitude of
possible such labels means that they are \textquotedblleft artificial
variables...added to the target states\textquotedblright\ (\cite%
{JoseTAES2018}, p. 1884). \ Rather, they are stand-ins for the unique
identities of actual physical objects and, as such, are no more
\textquotedblleft artificial\textquotedblright\ than any other signifier for
those objects. \ To claim otherwise is to implicitly claim that (for
example) target-classifier algorithms---which estimate the types or even
identities of physical targets---are \textquotedblleft
non-physical.\textquotedblright\ 

\begin{remark}
\label{Rem-Variable-ID}\textquotedblleft Target types\textquotedblright\ or
\textquotedblleft target classes\textquotedblright\ are signifiers that are
less accurate than identities but typically more accurate than labels. \ As
noted in \cite{Mah-Newbook}, p. 259, Remark 33, they---unlike
identities---are not necessarily time-invariant because some targets can
have multiple, distinct phenomenological \textquotedblleft
modes.\textquotedblright\ \ Typical examples include:\ variable swept-wing
aircraft (extended-wing vs. delta-wing); mobile missile launchers
(launch-ready vs. transiting); and diesel-electric submarines (surfaced vs.
submerged). \ This complication can be addressed within the RFS/LRFS\
framework but will not be further addressed in this paper.
\end{remark}

These points having been made, we are now in a position to address the
theoretically and phenomenologically correct mathematical representation of
time-evolving multitarget populations.

\subsection{Mathematical Representation of Multitarget States \label%
{A-TargRep-AA-MultiTarg}}

As noted earlier, the state of a point target at a particular instant has
the mathematical form \ $\xi =(\mathbf{x},\ell )\in \mathfrak{X}=\mathfrak{X}%
_{0}\times \mathfrak{L}$ \ where \ $\mathbf{x}\in \mathfrak{X}_{0}$ \ is its
kinematic state (position, velocity, orientation, etc.) and \ $\ell \in 
\mathfrak{L}$ \ is a provisional identifying label drawn (without
replacement) from a countable set \ $\mathfrak{L}$ \ of such labels. \ 

This, in turn, means that it is impossible for a target population to
contain pairs \ $(\mathbf{x},\ell )$, \ $(\mathbf{y},\ell )$ \ with \ $%
\mathbf{x}\neq \mathbf{y}$ \ since, then (for example) \textquotedblleft
Bob\textquotedblright\ could be in two different places simultaneously. \
Nor can it contain pairs \ $(\mathbf{x},\ell )$, \ $(\mathbf{x},\ell )$, \
since then two copies of (for example) \textquotedblleft
Bob\textquotedblright\ could exist simultaneously. \ Consequently, the state
of a target population must must be a collection \ $\xi _{1}=(\mathbf{x}%
_{1},\ell _{1}),...,\xi _{n}=(\mathbf{x}_{n},\ell _{n})$ \ where the
kinematic states have distinct labels---i.e., 
\begin{equation}
\ |\{(\mathbf{x}_{1},\ell _{1}),...,(\mathbf{x}_{n},\ell _{n})\}|=|\{\ell
_{1},...,\ell _{n}\}|=n.
\end{equation}

This collection does not have a natural ordering---e.g., is
\textquotedblleft Bob\textquotedblright\ greater or less than
\textquotedblleft Sue\textquotedblright ? \ Thus \ $\xi _{1},...,\xi _{n}$ \
must have the mathematical form \ $\{\xi _{1},...,\xi _{n}\}\subseteq 
\mathfrak{X}_{0}\times \mathfrak{L}$---i.e., a finite set---and not \ $(\xi
_{1},...,\xi _{n})\in (\mathfrak{X}_{0}\times \mathfrak{L})^{n}$---i.e., a
vector.

There are three additional reasons why vectors are theoretically
inappropriate representations of multitarget populations (\cite%
{MahSensors2019}, Section 2.4). \ 

\begin{enumerate}
\item \textit{Statistical Bias}. \ Imposing nonphysical information on a
physical system can create a statistical bias. \ The most obvious such
biases are deliberate---e.g., forcing a multitarget tracker-classifier to
prioritize targets of interest (ToI's), e.g., \textquotedblleft
Bob\textquotedblright\ is tactically more threatening than \textquotedblleft
Sue.\textquotedblright\ \ The concept of ToI is subjective and contextual,
not physical. \ Forcing it upon a target population results in a nonphysical
statistical bias in favor of more tactically significant ToI's. \ Similarly,
physical targets do not have a natural ordering. \ Forcing one upon a
population (as with vector representation)\ runs the risk of introducing an
unknown statistical bias.\ \ (In the TRFS\ framework, target ordering is
also viewed as \textquotedblleft non-physical,\textquotedblright\ see \cite%
{Svensson-FUSION2014}, Abstract.)

\item \textit{Uniqueness of States}. \ There should be a one-to-one
correspondence between physical states and their mathematical
representations. \ Because physical states have no inherent ordering, there
are \ $n!$ \ possible vector representations \ $\vec{\xi}_{\pi }=(\xi _{\pi
1},...,\xi _{\pi n})$ \ of the same multitarget population \ $\xi
_{1},...,\xi _{n}$, \ for permutations \ $\pi $ \ on \ $1,...,n$. \ There is
no way to choose one particular permutation as a representative, without
implicitly assuming that the targets have a specific ordering.

\item \textit{Performance Evaluation}. \ Multitarget tracking algorithms are
multitarget state estimators. \ Performance evaluation of such algorithms
requires the existence of a \textit{mathematical metric on multitarget states%
}, in order to measure the distance between the ground truth state and any
given estimated state. \ However, no such distance metric exists for vector
representation. \ Assume the contrary: \ a metric \ $d(\vec{\xi}_{1},\vec{\xi%
}_{2})$. \ Then if \ $\pi $ \ is not the identity permutation, \ $\vec{\xi}%
_{\pi }\neq \vec{\xi}$ \ and yet \ \ $d(\vec{\xi}_{\pi },\vec{\xi})=d(\vec{%
\xi},\vec{\xi})=0$, which contradicts the definition of a mathematical
metric on multitarget states.
\end{enumerate}

\begin{remark}
\label{Rem-GLMB}In objection to Item 1, it could be argued that the GLMB\
filter's labeling scheme is vulnerable to possible statistical bias. \ This
is not the case. \ As will be noted in Section \ref{A-TrajRep-AA-LRFS}, this
scheme assigns labels of the form \ $\ell _{1}=(k,1),...,\ell
_{n_{k}}=(k,n_{k})$ \ where \ $t_{k}$ \ is the time that an ensemble of \ $%
n_{k}$ \ targets appeared and where the natural numbers\ \ $1,2,...,n_{k}$ \
distinguish those\ targets from each other. \ But here is the crucial point:
\ It is not assumed that the list \ $1,...,n$ \ (and therefore also \ $\ell
_{1},...,\ell _{n_{k}}$) is ordered. \ The theoretical foundation of modern
mathematics is Cantor's set theory, whose primitive concept is the set. \
Vectors must therefore be defined in terms of sets. \ The simplest such
definition is as follows: \ $(x_{1},...,x_{n})$ \ is a finite set \ $%
\{(x_{1},1),...,(x_{n},n)\}$ \ that has been endowed with the ordering
relation \ $(x_{i},i)<(x_{j},j)$ \ if and only if \ $i<j$. \ If no such
relation has been stipulated (as in the GLMB scheme) then the list \ $%
1,...,n $ \ (and therefore also \ $\ell _{1},...,\ell _{n_{k}}$) has no
inherent ordering. \ Indeed, Greek letters, or any other arbitrary symbols
in any order, could have been used instead. \ A vector representation, on
the other hand, is not possible unless precisely this ordering relation has
been explicitly imposed on the collection \ $x_{1},...,x_{n}$.
\end{remark}

\section{Mathematical Representation of Trajectories \label{A-TrajRep}}

As noted in the Introduction, two mathematical representations of evolving
multitarget populations have been proposed: \ labeled\ RFS's (LRFS's) and
trajectory RFS's (TRFS's). \ The purpose of this section is to summarize and
contrast the two. \ It is organized as follows: \ LRFS's (Section \ref%
{A-TrajRep-AA-LRFS}); TRFS's (Section \ref{A-TrajRep-AA-TRFS}); comparison
of LRFS's and TRFS's (Section \ref{A-TrajRep-AA-Compare}); counterexamples
(Section \ref{A-TrajRep-AA-CounterEx}); and the trajectory PHD/CPHD\ filters
(Section \ref{A-TrajRep-AA-TPHD}).

\subsection{Labeled RFS's \label{A-TrajRep-AA-LRFS}}

From its inception in 1997, the RFS approach has included uniquely
identifying target identities or labels as target state variables---see \cite%
{GMN}, pp.135,196-197 and \cite{Mah-Artech}, pp. 505-507. \ Because of
computational considerations, however, the first implementations of RFS
filters mostly did not address track labeling; and, when they did, employed
computationally expensive techniques such as track-to-track association. \
Later implementations (for example, those based on Gaussian mixture or
particle-system approximation) addressed trajectories in a computationally
tractable manner via heuristic label-propagation schemes---see, e.g., \cite%
{Mah-Newbook}, p. 244-250. \ 

The labeled RFS (LRFS) theory of B.-T. Vo and B.-N. Vo \cite%
{VoVoTSPconjugate}, introduced in 2011 \cite{Vo-ISSNIP12-Conjugate}, is the
first systematic, theoretically rigorous formulation of true multitarget
tracking. \ Their provably Bayes-optimal GLMB filter is the currently most
sophisticated LRFS tracking algorithm. \ As previously noted, its latest
implementations can simultaneously track over one million targets in
significant clutter in real time using off-the-shelf computing equipment 
\cite{Beard2020}. \ 

This section summarizes LRFS theory, following the discussion in Chapter 15
of \cite{Mah-Newbook}.

As in \cite{GMN}, pp.135,196-197, the state of a point target is \ $(\mathbf{%
x},\ell )\in \mathfrak{X}=\mathfrak{X}_{0}\times \mathfrak{L}$ \ where \ $%
\mathfrak{X}_{0}$ \ is the kinematic state space and \ $\mathfrak{L}$ \ is a
countable space of labels \ $\ell $. \ As was noted in Remark \ref{Rem-GLMB}%
, GLMB\ filter labels have the form \ $(k,i)$ \ where \ $k\geq 0$ \ denotes
the time \ $t_{k}$ \ that the target first appeared and where \ $i\geq 1$ \
distinguishes it from all \ others born at the same time. \ As noted in \cite%
{VoVoSPIE13}, labels\ can be extended to contain target-type or
target-identity information, thus permitting simultaneous target tracking
and classification/identification.

A finite subset \ $X=\{(\mathbf{x}_{1},\ell _{1}),...,(\mathbf{x}_{n},\ell
_{n})\}\subseteq \mathfrak{X}$ \ is said to be \textit{labeled} if \ $|X_{%
\mathfrak{L}}|=|X|$ \ where \ $X_{\mathfrak{L}}=\{\ell _{1},...,\ell _{n}\}$
\ is its set of labels. \ That is: \ every target is assumed to have a
unique identifying label. \ If all instantiations of an RFS \ $\Xi \subseteq 
\mathfrak{X}$ \ are labeled, then \ $\Xi $ \ is an LRFS.

\begin{remark}
\label{Rem-Convention}Note that labeled sets are allowed to contain pairs of
the form \ $(\mathbf{x},\ell _{1})$, $(\mathbf{x},\ell _{2})$ \ with \ $\ell
_{1}\neq \ell _{2}$. \ This is not a conceptual f{}law. \ It is physically
possible, for example, for two or more targets traveling in formation to be
so close together that, within the resolution of the sensor, they appear to
have identical kinematic states over an extended duration of time. \
Additionally, because point targets are mathematical idealizations that do
not have a physical extent, it is possible for two different point targets
to (for example) have the same position.
\end{remark}

Bayes-optimal multitarget tracking is accomplished via the \textit{labeled
multitarget recursive Bayes filter}%
\begin{equation*}
...\rightarrow f_{k-1}(X|Z_{1:k-1})\rightarrow f_{k}(X|Z_{1:k-1})\rightarrow
f_{k}(X|Z_{1:k})\rightarrow ...
\end{equation*}%
where \ $Z_{1:k}:Z_{1},...,Z_{k},...$ \ is the time-sequence of
measurement-sets collected by a single sensor at times \ $%
t_{1},...,t_{k},... $ and \ $f_{k}(X|Z_{1:k})$ \ is a probability
distribution on \ $X\subseteq \mathfrak{X}_{0}\times \mathfrak{L}$. \ It
must be the case that \ $f_{k}(X|Z_{1:k})=0$ \ if \ $X$ \ is not
labeled---i.e., if \ $X$ \ is physically impossible. \ Stated differently, \ 
$f_{k}(X|Z_{1:k})$ \ is an \textquotedblleft LRFS
distribution.\textquotedblright\ \ Since targets evolve in four-dimensional
space-time, their states have the form \ $(\mathbf{x},\ell ,t)\in \mathfrak{X%
}=\mathfrak{X}_{0}\times \mathfrak{L}\times \mathbb{R}^{+}$ \ where \ $t$ \
is a known constant. \ Since time is assumed to belong to a discrete
sequence \ $t_{1},...,t_{k},...$, abbreviate \ $(\mathbf{x},\ell ,t_{k})$ \
as \ \ $(\mathbf{x},\ell ,k)\in \mathfrak{X}_{0}\times \mathfrak{L}\times 
\mathbb{N}$.

The \textit{labeled multi-Bernoulli }(LMB)\ \textit{distribution} is a
simple example of an LRFS distribution. \ It has the form%
\begin{eqnarray}
&&f_{J}(X) \\
&=&\delta _{|X|,|X_{\mathfrak{L}}|}\left( \prod_{\ell \in J-X_{\mathfrak{L}%
}}(1-q_{\ell })\right) \left( \prod_{(\mathbf{x,}\ell )\in X}\mathbf{1}%
_{J}(\ell )\cdot q_{\ell }s_{\ell }(\mathbf{x})\right)  \notag
\end{eqnarray}%
where \ $J$ \ is a finite subset of \ $\mathfrak{L}$; \ and where \ $q_{\ell
}$ \ and \ $s_{\ell }(x)$ \ are, respectively, the existence probability and
spatial distribution of the target with label \ $\ell \in J$. \ Note that \ $%
f_{J}(X)=0$ \ unless \ $|X|=|X_{\mathfrak{L}}|$ \ and \ $X_{\mathfrak{L}%
}\subseteq J$.

At time \ $t_{k}$, the multitarget state-set is estimated from \ $%
f_{k}(X|Z_{1:k})$\ \ using a Bayes-optimal multitarget state estimator,
e.g., the joint multitarget (JoM) or marginal multitarget (MaM)
estimator---see \cite{Mah-Artech}. \ If \ 
\begin{equation*}
\hat{X}^{k}=\{(\mathbf{\hat{x}}_{k}^{1},\hat{\ell}_{k}^{1},k),...,(\mathbf{%
\hat{x}}_{k}^{\hat{n}_{k}},\hat{\ell}_{k}^{\hat{n}_{k}},k)\}
\end{equation*}%
with \ $|\hat{X}^{k}|=\hat{n}_{k}$ \ then, at time \ $t_{k}$, \ $\hat{n}_{k}$
\ is the estimated number of targets and \ $\hat{\ell}_{k}^{1},...,\hat{\ell}%
_{k}^{\hat{n}_{k}}$ \ are their estimated labels. \ If \ $\hat{X}^{1},...,%
\hat{X}^{k}$\ \ is the time-sequence of multitarget state-estimates from
time \ $t_{1}$ \ to time \ $t_{k}$, define \ 
\begin{equation}
\hat{X}_{\ell }^{i}=\left\{ 
\begin{array}{ccc}
\{(\mathbf{x},\ell ,i)\} & \text{if} & (\mathbf{x},\ell ,i)\in \hat{X}^{i}%
\text{ \ for some \ }\mathbf{x}\in \mathfrak{X}_{0} \\ 
\emptyset & \text{if} & \text{otherwise}%
\end{array}%
\right. .
\end{equation}%
The\ $\ell $-trajectory (a.k.a. \ $\ell {}$-track)\ is the time-sequence \ $%
\hat{X}_{\ell }^{1}...,\hat{X}_{\ell }^{k},...$ \ of singleton or empty sets
(which thereby accounts for the acquisition, dropouts, and reacquisitions of
the target with label \ $\ell $). \ The time-consecutive nonempty
subsequences of an $\ell $-trajectory are its \textit{track segments}. \ 

It follows that, in the LRFS framework, a track segment can be represented
as a vector of the form%
\begin{equation}
\left( (\mathbf{x}^{1},\ell ,k),(\mathbf{x}^{2},\ell ,k+1),...,(\mathbf{x}%
^{i},\ell ,k+i-1)\right)  \notag
\end{equation}%
where \ $t_{k}$ \ is the segment's initial time, \ $i\geq 1$ \ is its
length, and \ $\mathbf{x}^{1},...,\mathbf{x}^{i}$ \ are its kinematic states
at times \ $t_{k},...,t_{k+i-1}$.

Now let \textquotedblleft $A\approx B$\textquotedblright\ abbreviate the
phrase \textquotedblleft $A$ is notationally equivalent to $B$%
\textquotedblright\ \ (in the sense that $\ A$ \ and \ $B$ \ are
characterized by the same parameters). \ Then we can successively re-notate
the track segment as follows:\ 

\begin{eqnarray*}
&&\left( (\mathbf{x}^{1},\ell ,k),(\mathbf{x}^{2},\ell ,k+1),...,(\mathbf{x}%
^{i},\ell ,k+i-1)\right) \\
&\approx &(\ell ,(\mathbf{x}^{1},k),(\mathbf{x}^{2},k+1),...,(\mathbf{x}%
^{i},k+i-1)) \\
&\approx &(\ell ,k,\mathbf{x}^{1},\mathbf{x}^{2},...,\mathbf{x}^{i}) \\
&\approx &(\ell ,k,\mathbf{x}^{1:i}).
\end{eqnarray*}%
That is: \ the track segment can be equivalently notated as \ $T_{\ell
}=(\ell ,k,\mathbf{x}^{1:i})$ \ where \ $0\leq k\leq k_{\max }$ \ and \ $%
1\leq i\leq k_{\max }-k+1$; and where \ $t_{k_{\max }}$ \ is the end-time of
the scenario. \ \ 

\subsection{Trajectory Random Finite Sets \label{A-TrajRep-AA-TRFS}}

The TRFS framework was introduced in 2014 in \cite{Svensson-FUSION2014} and
subsequently elaborated in \cite{AngelFUSION18}, \cite%
{GarciaSvensson-TSP2019}, \cite{GarciaSvensson-arXiv2019}, \cite%
{GranstromFUSION18}, and \cite{SvenssenFUSION2018}. \ There it was claimed
that LRFS labels are deficient because:

\begin{enumerate}
\item they are \textquotedblleft non-physical\textquotedblright\ (\cite%
{Svensson-FUSION2014}, Abstract);

\item they have "...no direct and unambiguous connection to the physical
phenomena under observation...\textquotedblright\ (\cite{Svensson-FUSION2014}%
, p. 2, column 2);

\item they \textquotedblleft ...do...not represent an underlying physical
reality...\textquotedblright\ (\cite{Svensson-FUSION2014}, p. 3, top of
column 1); and

\item a \textquotedblleft ...multitude of labeling [sic] can be
developed...[but] there is no general way of distinguishing the merits of a
given scheme...\textquotedblright ;\ which \textquotedblleft
...contravene[s]...the usual notion of there being a `one-to-one
correspondence between physical states and their mathematical
representations'...\textquotedblright\ (\cite{Svensson-FUSION2014}, p. 2,
column 2), here citing p. 405 of \cite{Mah-Artech}.
\end{enumerate}

In addition, it was claimed in \cite{Svensson-FUSION2014} that the example
in Figure 2 of four slightly different scenarios with \textquotedblleft
...two targets approaching on the real line, pausing and then separating
illustrates the appeal of using RFS trajectories [rather than
labels]...\textquotedblright\ (\cite{Svensson-FUSION2014}, column 1). \ 

\textquotedblleft Trajectory random finite sets\textquotedblright\ (TRFS's)
were therefore proposed as a replacement for LRFS's. \ As described in pp.
1-2 of \cite{AngelFUSION18}, a \textquotedblleft
trajectory\textquotedblright\ is what we earlier called a track segment. \
It has the form \ $T=(k,\mathbf{x}^{1:k})$ \ where \ $0\leq k\leq k_{\max }$
\ is the trajectory's initial time, \ $1\leq i\leq k_{\max }-k+1$ \ is its
length, and \ $\mathbf{x}^{1},...,\mathbf{x}^{i}$ \ are its kinematic states
at times \ $t_{k},...t_{k+i-1}$. \ A finite set of trajectories (SoT) is,
therefore, \ $\mathbf{T}=\{T_{1},...,T_{n}\}$ \ where \ $T_{1},...,T_{n}$ \
are trajectories. \ A TRFS is a random variable on the class of SoT's. \ 

It was also claimed that the TRFS approach subsumes the LRFS approach as a
special case because \ $\mathbf{T}=\{(k,\mathbf{x}_{1}),...,(k,\mathbf{x}%
_{n})\}$ \ is an equivalent representation of a set \ $X_{k}=\{\mathbf{x}%
_{1},...,\mathbf{x}_{n}\}$ \ of targets at time \ $t_{k}$.\ 

The TRFS framework subsequently provided the basis for the trajectory PHD\
(TPHD)\ filter\ \cite{AngelFUSION18}; the trajectory CPHD (TCPHD)\ filter\ 
\cite{GarciaSvensson-TSP2019}; and the \textquotedblleft trajectory Poisson
multi-Bernoulli mixture trajectory filter\textquotedblright\ \cite%
{GranstromFUSION18}, \cite{SvenssenFUSION2018}. \ The trajectory PHD and
CPHD\ filters will be discussed in detail in Section \ref{A-TrajRep-AA-TPHD}.

\subsection{Comparison of the LRFS and TRFS Frameworks \label%
{A-TrajRep-AA-Compare}}

The four enumerated criticisms of the LRFS framework in the previous section
have already been rebutted in Section \ref{A-TargRep-AA-Ident}. \ To
reiterate: \ Neither the provisionality nor the infinitude of track labels
means that they are \textquotedblleft non-physical.\textquotedblright\ \
Rather, they are stand-ins for the unique identities of actual physical
objects. \ As such, they are no more non-physical\ than any other signifier
for those objects.

As for the argument relating to Figure 2 of \cite{Svensson-FUSION2014}, it
appears to involve an apples-with-oranges comparison. \ It was implicitly
assumed there that the four trajectories were all being considered as
entireties---that is, the four scenarios ended before measurement processing
commenced. \ However, it is well known that a smoother algorithm will
typically outperform a filter algorithm, because the former processes both
past and future measurements whereas the latter can process past
measurements only. \ Stated differently, the smoother is a batch-processing
algorithm whereas the filter cannot discern the future and thus is, in this
sense, \textquotedblleft real-time.\textquotedblright\ \ TRFS algorithms,
such as the TPHD and TCPHD filters discussed in Section \ref%
{A-TrajRep-AA-TPHD}, are also of the smoother type since, at any given time
\ $t_{k}$, all trajectories prior to \ $t_{k}$ \ must be taken into
consideration.

It might thereby be argued that the TRFS framework is still valid within the
context of multitarget smoothing. \ But as we will shortly see, this is not
the case.

Let us now turn to a more detailed comparison of the LRFS and TRFS
frameworks. \ In the former, a \textquotedblleft
trajectory\textquotedblright\ (track segment) has the form \ $T_{\ell
}=(\ell ,k,\mathbf{x}^{1:i})$; \ whereas in the latter it has the form \ $%
T=(k,\mathbf{x}^{1:i})$. \ Since a SoT has the form \ $\mathbf{T}%
=\{T_{1},...,T_{n}\}$, \ it must therefore be the case that \ $T_{j}=(k_{j},%
\mathbf{x}_{j}^{1:i_{j}})$ \ for some \ $k_{j}$, \ $i_{j}$, \ and \ $\mathbf{%
x}_{j}^{1},...,\mathbf{x}_{j}^{i_{j}}$. \ 

That is: \ \textit{the trajectory} \ $T_{j}$ \ \textit{has been implicitly
assigned the integer label} \ $j$---\textit{which is then ignored
(\textquotedblleft stripped off\textquotedblright )}. \ But as was noted in
Section \ref{A-TargRep-AA-Ident}, when the signifier of an entity is
jettisoned, critical information about it is jettisoned as well. \ As the
following subsection demonstrates, this loss of information results in
numerous and serious difficulties.

\subsection{Counterexamples \label{A-TrajRep-AA-CounterEx}}

\textit{Counterexample 1}: \ The LRFS framework is not a special case of the
TRFS framework. \ Let \ $\mathbf{T}=\{(k,\mathbf{x}_{1}),...,(k,\mathbf{x}%
_{n})\}$ \ with \ $\mathbf{x}_{1},...,\mathbf{x}_{n}$ \ distinct. \ Then
this is mathematically equivalent to the LRFS\ representation \ $%
X_{k}=\{(1,k,\mathbf{x}_{1}),...,(n,k,\mathbf{x}_{n})\}$. \ From this one
might conclude that the TRFS framework includes the LRFS\ framework as a
special case. \ But this is not so. \ In Remark \ref{Rem-Convention} it was
noted that a labeled set of the form \ $X=\{(\mathbf{x},\ell _{1}),(\mathbf{x%
},\ell _{2})\}$ \ with \ $\ell _{1}\neq \ell _{2}$ \ is possible. \ Thus
consider the case \ $\mathbf{x}_{1}=...=\mathbf{x}_{n}=\mathbf{x}$, \ in
which case \ $\mathbf{T}=\{(k,\mathbf{x})\}$ \ whereas \ $X_{k}=\{(1,k,%
\mathbf{x}),...,(n,k,\mathbf{x})\}$ \ is a valid labeled finite set. \ That
is: \ there is no way to represent such labeled multitarget states in the
TRFS framework. \ And this is not the only possible such anomaly---see
Counterexample 3.

\textit{Counterexample 2}: \ Contrary to claim, the TRFS framework
contravenes the usual notion of there being a one-to-one correspondence
between physical states and their mathematical representations. \ Consider \ 
$\mathbf{T}_{0}=\{T_{0}\}$ \ and \ $\mathbf{T}_{1}=\{T_{1},T_{2},T_{3}\}$ \
where \ $T_{0}=(k,\mathbf{x},\mathbf{x}^{1},\mathbf{x}^{2})$, \ $T_{1}=(k,%
\mathbf{x})$, $T_{2}=(k+1,\mathbf{x}^{1})$, \ $T_{3}=(k+2,\mathbf{x}^{2})$.
\ Then because \ $(k,\mathbf{x},\mathbf{x}^{1},\mathbf{x}^{2})$ \ is an
abbreviation of \ $((k,\mathbf{x}),(k+1,\mathbf{x}^{1}),(k+2,\mathbf{x}%
^{2})) $, \ it follows that \ $\mathbf{T}_{0}$ \ and \ $\mathbf{T}_{1}$ \
are mathematically distinct representations of the same physical trajectory
\ $(k,\mathbf{x})$, $(k+1,\mathbf{x}^{1})$, $(k+2,\mathbf{x}^{2})$. \ Now
restore the implicit labels that have been stripped off: \ $T_{0}=(0,k,%
\mathbf{x},\mathbf{x}^{1},\mathbf{x})$, \ $T_{1}=(1,k,\mathbf{x})$, $\
T_{2}=(2,k+1,\mathbf{x}^{1})$, \ $T_{3}=(3,k+2,\mathbf{x}^{2})$. \ Then the
difficulty vanishes because \ $\mathbf{T}_{0}$ \ represents the trajectory
of a single physical target with label \ $0$; whereas \ $\mathbf{T}_{1}$ \
represents the trajectories of three successively appearing and disappearing
targets with respective labels \ $1$, $2$, $3$.

\textit{Counterexample 3}: \ Impossible scenarios can be represented in the
TRFS framework even while valid ones cannot. \ Consider \ $\mathbf{T}%
_{2}=\{T_{1},T_{2}\}$ \ where $T_{1}=(k,\mathbf{x},\mathbf{x}^{1})$ \ and $%
T_{2}=(k,\mathbf{x},\mathbf{x}^{2})$\ \ with \ $\mathbf{x},\mathbf{x}^{1},%
\mathbf{x}^{2}$ \ distinct. \ Then \ $\mathbf{T}_{2}$ \ is physically
impossible since a single target \ $\mathbf{x}$ \ at time \ $t_{k}$ \ cannot
evolve to two different states \ $\mathbf{x}^{1}$ and \ $\mathbf{x}^{2}$ \
at time \ $t_{k+1}$. \ Now restore the stripped labels: \ $T_{1}=(1,k,%
\mathbf{x},\mathbf{x}^{1})$, $\ T_{2}=(2,k,\mathbf{x},\mathbf{x}^{2})$. \
Then \ $\mathbf{T}_{2}$ \ represents a target separation (e.g., a
target-spawning event): \ targets \ $1,2$ \ had identical states \ $\mathbf{x%
}$ \ at time \ $t_{k}$, at which point they separated and evolved
respectively to \ $\mathbf{x}^{1}$ \ and \ $\mathbf{x}^{2}$ \ at time \ $%
t_{k+1}$. \ 

\textit{Counterexample 4}: \ Contrary to claim, stripping labels from states
increases tracking uncertainty. \ Consider \ $\mathbf{T}_{3}=\{T_{1},T_{2}\}$
\ where \ $T_{1}=(k,\mathbf{x}^{1:5})$ \ and \ $T_{2}=(k+10,\mathbf{y}%
^{1:5}) $. \ Then there is an ambiguity: \ Does \ $\mathbf{T}_{3}$ \
represent a single dropped and then reacquired track, or two successively
appearing and disappearing tracks? \ This ambiguity is resolved if we
restore stripped labels: \ either \ $T_{1}=(1,k,\mathbf{x}^{1:5})$ \ and \ $%
T_{2}=(1,k+10,\mathbf{x}^{1:5})$ \ (if restored to a single reacquired
track) or \ $T_{1}=(1,k,\mathbf{x}^{1:5})$ \ and \ $T_{2}=(2,k+10,\mathbf{y}%
^{1:5})$ \ (if restored to two consecutive tracks).\ 

\subsection{Trajectory PHD and CPHD\ Filters \label{A-TrajRep-AA-TPHD}\ \ }

These are not the only difficulties with the TRFS framework. \ As noted
earlier, the TPHD/TCPHD filters were introduced in \cite{AngelFUSION18}\ and 
\cite{GarciaSvensson-arXiv2019}, \cite{GarciaSvensson-TSP2019},
respectively. \ The purpose of this subsection is to demonstrate that they
are theoretically problematic. \ It is organized as follows: \ C/PHD\
filters (Section \ref{A-TrajRep-AA-TPHD-AAA-PHD}); multitarget state
estimation for C/PHD\ filters (Section \ref{A-TrajRep-AA-TPHD-AAA-E}); the
TPHD\ filter (Section \ref{A-TrajRep-AA-TPHD-AAA-TPHD}); multitarget state
estimation for the TPHD\ filter (Section \ref{A-TrajRep-AA-TPHD-AAA-TE});
and the theoretical basis of the TPHD filter (Section \ref%
{A-TrajRep-AA-TPHD-AAA-Basis}).

\subsubsection{The C/PHD\ Filters \label{A-TrajRep-AA-TPHD-AAA-PHD}}

The conventional PHD filter \cite{Mah-Artech} has the form \ 
\begin{equation*}
...\rightarrow D_{k-1}(\mathbf{x}|Z_{1:k-1})\rightarrow D_{k}(\mathbf{x}%
|Z_{1:k-1})\rightarrow D_{k}(\mathbf{x}|Z_{1:k})\rightarrow ...
\end{equation*}%
Here,\ \ $D_{k}(\mathbf{x}|Z_{1:k})$ \ with \ $\mathbf{x}\in \mathfrak{X}%
_{0} $ \ is the first-order statistical moment (a.k.a. PHD or
\textquotedblleft intensity density\textquotedblright ) of \ $%
f_{k}(X|Z_{1:k})$ \ with \ $X\subseteq \mathfrak{X}_{0}$, \ which in turn is
assumed to be the distribution of a Poisson RFS with PHD \ $D_{k}(\mathbf{x}%
|Z_{1:k})$. \ 

The CPHD\ filter \cite{Mah-Artech}\ propagates the cardinality distribution
\ $p_{k}(n|Z_{1:k})$ \ (i.e., the probability that there are \ $n$ \ targets
present at time \ $t_{k}$) in addition to the PHD \ $D_{k}(\mathbf{x}%
|Z_{1:k})$. \ It is based on i.i.d.c. RFS's (which are generalizations of
Poisson\ RFS's). \ 

\begin{remark}
\label{Rem-Poisson}Poisson RFS's on \ $\mathfrak{X}_{0}\times \mathfrak{L}$
\ are not LRFS's because their distributions have the form%
\begin{equation}
f(X)=e^{-\sum_{\ell }\int D(\mathbf{x},\ell )d\mathbf{x}}\prod_{(\mathbf{x}%
,\ell )\in X}D(\mathbf{x},\ell ).
\end{equation}%
Suppose that \ $Y=\{(\mathbf{x}_{1},\ell ),(\mathbf{x}_{2},\ell )\}$ \ for
some \ $\mathbf{x}_{1}\neq \mathbf{x}_{2}$ \ such that $\ D(\mathbf{x}%
_{1},\ell )>0$, \ $D(\mathbf{x}_{2},\ell )>0$. \ Then \ $Y$ \ is physically
impossible since \ $|Y|=2\neq 1=|Y_{\mathfrak{L}}|$. \ Yet \ $f(Y)\propto D(%
\mathbf{x}_{1},\ell )\cdot D(\mathbf{x}_{2},\ell )\neq 0$---i.e., \ $f(X)$ \
is not an LRFS distribution.
\end{remark}

\subsubsection{Multitarget State Estimation for the PHD\ Filter \label%
{A-TrajRep-AA-TPHD-AAA-E}}

This proceeds as follows. \ First, compute 
\begin{equation}
\bar{N}_{k}=\int D_{k}(\mathbf{x}|Z_{1:k})d\mathbf{x,}
\end{equation}%
which is the expected number of targets in the multitarget population at
time \ $t_{k}$. \ Second, round \ $\bar{N}_{k}$ \ off to the nearest integer
\ $\hat{n}_{k}$. \ Third, determine the states \ $\mathbf{\hat{x}}_{1},...,%
\mathbf{\hat{x}}_{\hat{n}_{k}}$ \ such that\ \ $D_{k}(\mathbf{\hat{x}}%
_{1}|Z_{1:k})$,...,$D_{k}(\mathbf{\hat{x}}_{\hat{n}_{k}}|Z_{1:k})$ \ are the
heights of the \ $\hat{n}_{k}$ \ tallest \textquotedblleft
peaks\textquotedblright\ of the graph of the function \ $D_{k}(\mathbf{x}%
|Z_{1:k})$.

State estimation for the CPHD\ filter differs only in that \ 
\begin{equation}
\hat{n}_{k|k}=\arg \sup_{n\geq 0}p_{k}(n|Z_{1:k}).
\end{equation}

\subsubsection{The TPHD\ Filter \label{A-TrajRep-AA-TPHD-AAA-TPHD}\ }

In \cite{AngelFUSION18} it was claimed that the PHD filter can be directly
generalized from RFS's to TRFS's. \ This \textquotedblleft trajectory PHD
filter\textquotedblright\ (TPHD filter) has the form%
\begin{equation*}
...\rightarrow \vec{D}_{k-1}(T|Z_{1:k-1})\rightarrow \vec{D}%
_{k}(T|Z_{1:k-1})\rightarrow \vec{D}_{k}(T|Z_{1:k})\rightarrow ...
\end{equation*}%
where (in the case of \ $\vec{D}_{k}(T|Z_{1:k})$) \ $T=(k^{\prime },\mathbf{x%
}^{1:i})$ \ for all \ $0\leq k^{\prime }\leq k$, \ $1\leq i\leq k-k^{\prime
}+1$, and \ $\mathbf{x}^{1:i}\in \mathfrak{X}^{i}$. \ This means that all
trajectories prior to time \ $t_{k}$---i.e., all $(k^{\prime },\mathbf{x}%
^{1:i})$ \ with \ $0\leq k^{\prime }\leq k-1$ \ and \ $1\leq i\leq
k-k^{\prime }$---must be taken into consideration. \ Consequently,
calculation of the TPHD \ $\vec{D}_{k}(T|Z_{1:k})$ \ for all \ $T$ \ is
significantly more computationally intensive than calculation of the
corresponding \ PHD \ $D_{k}(\mathbf{x}|Z_{1:k})$ \ for all \ $\mathbf{x}$.

\subsubsection{Multitarget State Estimation for the TPHD\ Filter \label%
{A-TrajRep-AA-TPHD-AAA-TE}}

This is claimed to have exactly the same form as before, i.e.: \ (a)\
compute (\cite{AngelFUSION18}, Eq. (1)) \ 
\begin{eqnarray}
\bar{N}_{k} &=&\int \vec{D}_{k}(T|Z_{1:k})dT \\
&=&\sum_{k^{\prime }=0}^{k}\sum_{i=1}^{k-k^{\prime }+1}\vec{D}_{k}(k^{\prime
},\mathbf{x}^{1:i}|Z_{1:k})d\mathbf{x}^{1:i};
\end{eqnarray}%
\ (b) round off \ $\bar{N}_{k}$ \ to \ $\hat{n}_{k}$; \ and (c)\ determine
the locations \ $\hat{T}_{1},...,\hat{T}_{\hat{n}_{k}}$ \ of the \ $\hat{n}%
_{k}$ \ largest peaks of the graph of \ $\vec{D}_{k}(T|Z_{1:k})$. \ 

\textit{This TPHD/TCPHD estimation procedure is mathematically undefined}.\
\ For, let \ $\iota $ \ be the unit of measurement of the single-target
state space \ $\mathfrak{X}$. \ Then the unit of measurement of a trajectory
\ $(k^{\prime },\mathbf{x}^{1:i})$ \ is \ $\iota ^{i}$ \ and so the unit of
\ $\vec{D}_{k}(k^{\prime },\mathbf{x}^{1:i}|Z_{1:k})$ \ is \ $\iota ^{-i}$.
\ Because it is impossible to numerically compare \ $\vec{D}_{k}(k^{\prime },%
\mathbf{x}^{1:i}|Z_{1:k})$ \ to \ $\vec{D}_{k}(k^{\prime \prime },\mathbf{x}%
^{1:i^{\prime }}|Z_{1:k})$ \ when \ $i\neq i^{\prime }$, \ it is impossible
to determine \ $\hat{T}_{1},...,\hat{T}_{\hat{n}_{k}}$. \ 

This remains true when the TPHD/TPHD filters are implemented using Gaussian
mixture approximations. \ These have the form (\cite{AngelFUSION18}, Eqs.
(16,17)):%
\begin{equation}
\vec{D}_{k}(k^{\prime },\mathbf{x}^{1:i})\cong
\sum_{j=1}^{N_{k}}w_{k,j}\cdot \vec{N}_{\vec{P}_{k,j}^{1:i}}((k^{\prime },%
\mathbf{x}^{1:i});(k_{k,j}^{\prime },\mathbf{x}_{k,j}^{1:i_{k,j}}))
\label{eq-GM}
\end{equation}%
where \ $w_{k,j}\geq 0$; \ where%
\begin{equation}
\vec{N}_{\vec{P}_{0}^{1:i_{0}}}((k^{\prime },\mathbf{x}^{1:i});(k_{0},%
\mathbf{x}_{0}^{1:i_{0}}))=\delta _{k^{\prime },k_{0}}\delta _{i,i_{0}}\cdot
N_{\vec{P}_{0}^{1:i_{0}}}(\mathbf{x}^{1:i_{0}}-\mathbf{x}_{0}^{1:i_{0}})
\end{equation}%
is a probability distribution on trajectories \ $T=(k^{\prime },\mathbf{x}%
^{1:i})$; and where \ $N_{\vec{P}_{0}^{1:i}}(\mathbf{x}^{1:i})$ \ denotes a
zero-mean Gaussian distribution in the variable \ $\mathbf{x}^{1:i}=(\mathbf{%
x}^{1},...,\mathbf{x}^{i})$ \ with covariance matrix \ $\vec{P}_{0}^{1:i}$.
\ 

We may assume that \ $w_{k,j}>0$ \ for all \ $j$, which forces \ $i_{k,j}=i$
\ for all \ $j$ \ if the summation in (\ref{eq-GM}) is to be mathematically
well-defined. \ Thus at its most general, the GM in (\ref{eq-GM}) must have
the form 
\begin{equation}
\vec{D}_{k}(k^{\prime },\mathbf{x}^{1:i})\cong \bar{N}_{k}%
\sum_{j=1}^{N_{k}}w_{k,k^{\prime },i,j}\cdot \vec{N}_{\vec{P}_{k,j}^{1:i}}(%
\mathbf{x}^{1:i}-\mathbf{x}_{k,k^{\prime },i,j}^{1:i})
\end{equation}%
where \ $\bar{N}_{k}>0$ \ and where \ $w_{k,k^{\prime },i,j}>0$ \ are such
that \ 
\begin{equation}
\sum_{k^{\prime }=0}^{k}\sum_{i=1}^{k-k^{\prime
}+1}\sum_{j=1}^{N_{k}}w_{k,k^{\prime },i,j}=1,
\end{equation}%
from which follows \ $\int \vec{D}_{k}(T)dT\cong \bar{N}_{k}$. \ Since \ $%
\vec{D}_{k}(k^{\prime },\mathbf{x}^{1:i})$ \ and \ $\vec{D}_{k}(k^{\prime
\prime },\mathbf{x}^{1:i^{\prime }})$ \ are numerically incommensurable when
\ $i\neq i^{\prime }$, \ it follows that the inequality \ $w_{k,k^{\prime
},i,j}>w_{k,k^{\prime \prime },i^{\prime },j^{\prime }}$ \ cannot compel us
to conclude that GM component \ $\mathbf{x}_{k,k^{\prime },i,j}^{1:i}$ \ has
a higher \textquotedblleft peak\textquotedblright\ than GM component \ $%
\mathbf{x}_{k,k^{\prime \prime },i^{\prime },j^{\prime }}^{1:i^{\prime }}$.
\ That is: it is not true that \textquotedblleft ...the estimated set of
trajectories corresponds to...the components with highest
weights...\textquotedblright\ (as was asserted in regard to Eq. (21) of \cite%
{AngelFUSION18}).

\subsubsection{Theoretical Basis of the TPHD\ Filter \label%
{A-TrajRep-AA-TPHD-AAA-Basis}}

The TPHD filter is inherently theoretically erroneous because it requires
\textquotedblleft Poisson TRFS's\textquotedblright ---which, like Poisson
LRFS's, do not exist. \ For, according to (\cite{AngelFUSION18}, Eq. (7))
the distribution of such a TRFS has the form \ \ 
\begin{equation}
f(\mathbf{T})=e^{-\int \vec{D}(T)dT}\prod_{T\in \mathbf{T}}\vec{D}(T).
\end{equation}%
Let \ $\mathbf{T}_{01}=\mathbf{T}_{0}\uplus \mathbf{T}_{1}$ \ where \ $%
\mathbf{T}_{0},\mathbf{T}_{1}$ \ were defined in Counterexample 2 of Section %
\ref{A-TrajRep-AA-CounterEx}. \ Then \ $\mathbf{T}_{01}$ \ is a physically
impossible SoT since it contains two distinct instances of the same physical
trajectory. \ If \ $f(\mathbf{T})$ \ is a TRFS distribution then it must
vanish on impossible SoT's---in particular, it must be the case that \ $f(%
\mathbf{T}_{01})=0$. \ But 
\begin{equation}
f(\mathbf{T}_{01})\propto \vec{D}(T_{0})\cdot \vec{D}(T_{1})\cdot \vec{D}%
(T_{2})\cdot \vec{D}(T_{3})\neq 0,
\end{equation}%
where the four factors on the right can be nonzero since \ $%
T_{0},T_{1},T_{2},T_{3}$ \ are, individually, valid trajectories. \ 

The same comments apply to the trajectory Poisson multi-Bernoulli mixture
trajectory filter\ \cite{GranstromFUSION18}, \cite{SvenssenFUSION2018} since
it also requires non-existent Poisson TRFS's. \ The TCPHD\ filter is
similarly erroneous since it requires i.i.d.c. TRFS's---which do not exist
for the same reason that Poisson TRFS's do not exist.

\section{Conclusions \label{A-Concl} \ }

This paper described and compared two proposed mathematical representations
of multiple-target trajectories: \ the LRFS vs. TRFS frameworks. \ It was
shown that the latter is questionable in multiple and serious respects,
whereas the claimed deficiencies of the former are mistaken. \ This was
accomplished via a systematic, step-by-step analysis beginning with
phenomenological and epistemological basics---in particular, the semiotic
concepts of \textit{signified} vs. \textit{signifier}.

\end{document}